\documentclass[12pt]{article}
\usepackage{latexsym, amssymb, amsmath, amscd, amsfonts, epsfig, graphicx, colordvi,verbatim,ifpdf,tikz,extarrows}
\usepackage{amsfonts, amsmath, amssymb}
\usepackage{amssymb,amsfonts,amsmath,latexsym,epsfig,cite, psfrag,eepic,color}
\usepackage{amscd,graphics}
\usepackage{latexsym, amssymb,  amsmath,amscd, amsfonts, epsfig, graphicx, colordvi,amsthm}
\usepackage{ytableau}

\usepackage[colorlinks,linkcolor=blue,anchorcolor=blue,citecolor=blue]{hyperref}

\usepackage{graphicx}
\usepackage{color}
\usepackage{ifpdf}
\usepackage{fancybox}

\usepackage{float}

\newtheorem{thm}{Theorem}[section]

\newtheorem{lem}[thm]{Lemma}

\def\pf{\noindent{\it Proof.} }
\def\qed{\nopagebreak\hfill{\rule{4pt}{7pt}}
\medbreak}

\hypersetup{colorlinks = true}

\numberwithin{equation}{section}

\def\qed{\nopagebreak\hfill{\rule{4pt}{7pt}}
\medbreak}

\newlength{\boxedparwidth}
  {\begin{center} \begin{tabular}{|@{\hspace{.315in}}c@{\hspace{.15in}}|}
                  \hline \\ \begin{minipage}[t]{\boxedparwidth}
                  \setlength{\parindent}{.25in}}%
  {\end{minipage} \\ \\ \hline \end{tabular} \end{center}}

\begin{document}
\begin{center}

 {\Large \bf The inequality on the number of $1$-hooks, $2$-hooks and $3$-hooks in $t$-regular partitions}
\end{center}

\begin{center}
  {Hongshu Lin}, and {Wenston J.T. Zang$^2$} \vskip 2mm

 $^{1,2}$School of Mathematics and Statistics, Northwestern Polytechnical University, Xi'an 710072, P.R. China\\[6pt]
	$^{1,2}$ MOE Key Laboratory for Complexity Science in Aerospace, Northwestern Polytechnical University, Xi'an 710072, P.R. China\\[6pt]
	$^{1,2}$ Xi'an-Budapest Joint Research Center for Combinatorics, Northwestern Polytechnical University, Xi'an 710072, P.R. China\\[6pt]

   \vskip 2mm

$^1$linhongshu@mail.nwpu.edu.cn,
    $^2$zang@nwpu.edu.cn
\end{center}

\vskip 6mm \noindent {\bf Abstract.} Let $b_{n,k}$ denote the number of hooks of length $k$ in all the $t$-regular partitions of $n$. Singh and Barman raised the question of finding the relation between $b_{t,2}(n)$ and $b_{t,1}(n)$. Kim showed that there exists $N$ such that $b_{t,2}(n)\ge b_{t,1}(n)$ and $b_{t,2}(n) \geq b_{t,3}(n)$ for $n>N$. In this paper, we find an explicit bound of $N=O(t^5)$ for $b_{t,2}(n)\geq b_{t,1}(n)$ and show that $b_{t,2}(n) \geq b_{t,3}(n)$ for all $n\ge 4$.



\section{Introduction}


This paper focuses on the hook length biases in $t$-regular partitions. The partition of positive integers $n$ is a finite sequence of non-increasing positive integers $\{ \lambda_{1},\lambda_{2},\dots,\lambda_{\ell} \} $ such that $\lambda_{1} + \lambda_{2}+ \dots + \lambda_{\ell} = n $. A $t$-regular partition of a positive integer $n$ is a partition of $n$ such that none of its parts is divisible by $t$. We used $b_t(n)$ to denote the number of $t$-regular partitions of $n$.

\begin{figure}[h]
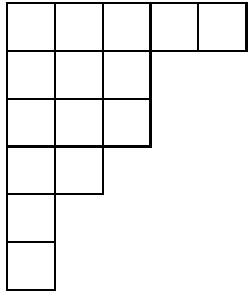

		\begin{center}
		\begin{ytableau}[]
        \ & \ & \ & \ & \ \\
        \ & \ & \ \\
        \ & \ & \ \\
        \ & \ \\
        \ \\
        \
        \end{ytableau}
		\caption{The Young diagram of partition $(5,3,3,2,1,1)$}
		\label{F1}
		\end{center}
	\end{figure}

\begin{figure}[h]
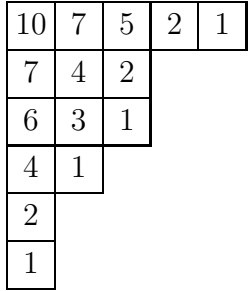

		\begin{center}
		\begin{ytableau}[]
        10 & 7 & 5 & 2 & 1\\
        7 & 4 & 2 \\
        6 & 3 & 1 \\
        4 & 1 \\
        2 \\
        1
        \end{ytableau}
		\caption{The hook lengths of $(5,3,3,2,1,1)$}
		\label{F2}
		\end{center}
	\end{figure}

A Young diagram of a partition  $\{ \lambda_{1},\lambda_{2},\dots,\lambda_{r} \} $ is a left-justified array of boxes with i-th row (from the top) having $\lambda_{i}$ boxes.The hook length of a box in a Young diagram is the sum of the number of the boxes directly right to it, the number of boxes directly below it and 1 (for the box itself). For example, the Young diagram of partition $(5,3,3,2,1,1)$ is shown in Figure \ref{F1}. The hook length of a box in the Young diagram is the number of the boxes directly to its right or directly below it and including itself exactly once. For example, Figure \ref{F2} shows the hook length of each box in the Young diagram of partition $(5,3,3,2,1,1)$.


Hook length plays an important role in the theory of representation theory of the symmetric group $S_n$ and the general linear group $GL_n(\mathbb{C})$. Hook length also relates to the theory of symmetric functions as the number of the standard Young tableaux, which is well known as the hook-length formula, see \cite{Frame-Rhobinson-Thrall-1954,Young-1927} for more details. The Nekrasov-Okounkov formula (see \cite{Han-Nek-2010, Nek-Oko}) builds the connection between hook length and the theory of modular forms. There are many other studies on the hook length, including the asymptotic, combinatorial and arithmetic properties, especially the $t$-core partitions and the $t$-hook partitions, see \cite{Garvann-Kim, James-Kerber, Littlewood-Modular} for example.

For integers $t \geq 2 $ and $ k \geq 1$, let $b_{t,k}(n)$ denote the number of hooks of length $k$ in all the $t$-regular partitions of $n$. Singh and Barman \cite{Singh-Barman-2024} proved that $b_{t+1,1}(n)\geq b_{t,1}(n)$ for all $t\geq 2 $ and $n\geq 0$. They also found that $b_{2,2}(n)\ge b_{2,1}(n)$ for all $n > 4$ and $b_{2,2}(n)\ge b_{2,3}(n)$ for all $n\ge0$. Moreover, they  proposed the conjecture of $b_{3,2}(n)\ge b_{3,1}(n)$ for $n\ge 28$. This conjecture was recently proved by Qu and Zang \cite{qu-zang-2025}.

At the end of their paper\cite{Singh-Barman-2024}, Singh and Barman mentioned that it would be interesting to study the relationship between $b_{t,1}(n)$ and $b_{t,2}(n)$. The asymptotic relationship between $b_{t,1}(n)$ and $b_{t,2}(n)$ was given by Kim \cite{Kim}. Kim  found  the asymptotic formula for  $b_{t,1}(n)$, $b_{t,2}(n)$ and $b_{t,3}(n)$, and   showed the following inequalities for sufficiently large $n$.

\begin{thm}[\cite{Kim}]
Let $t\geq 2$ be an integer. For sufficiently large integers $n$,
\begin{equation}
    b_{t,2}(n)\geq b_{t,1}(n) \text{ and } b_{t,2}(n) \geq b_{t,3}(n).
\end{equation}
\end{thm}
 The main  result of this paper is to find an explicit bound of the above theorem. To be specific, on the inequality $b_{t,2}(n)\geq b_{t,1}(n)$, we have the following result.

\begin{thm}\label{thm-com-int-umn}
  For integers $ n\geq 192t^5-192t^4-24t^3+24t^2+6t+2$, we have
  \begin{equation}\label{equ-thm-main}
 b_{t,2}(n) \geq b_{t,1}(n).
\end{equation}
\end{thm}

For the inequality $b_{t,2}(n)\geq b_{t,3}(n)$, we may give a necessary and sufficient condition as given below.
\begin{thm}\label{thm-bt2-bt3}
The inequality
\begin{equation}
    b_{t,2}(n) \geq b_{t,3}(n).
\end{equation}
holds if and only if the integer pair $(t,n)$ satisfies either $t=2$ or $n\ne 3$.
\end{thm}

This paper is organized as follows. In Section \ref{2}, we give a proof of Theorem \ref{thm-com-int-umn}. To this end, we divide the
difference between the generating function of $b_{t,2}(n) $ and $ b_{t,1}(n) $ into two parts and find the non-negativity of each part for specific $n$. We prove Theorem \ref{thm-bt2-bt3} in Section \ref{3} by dividing the
 generating function of $b_{t,2}(n)-b_{t,3}(n) $ into three parts. We then give combinatorial proofs that each part is non-negative for specific $n$.

\section{A proof of Theorem \ref{thm-com-int-umn}}\label{2}

In this section, we give a proof of Theorem \ref{thm-com-int-umn}.  We first write the  generating function of $b_{t,2}(n)-b_{t,1}(n) $ as a summation of two summands, which will be  stated in Theorem \ref{-A+B+C}. We then show that each summand is non-negative for specific $n$ through combinatorial method.

We first introduce some notations on integer partitions. Here and throughout this paper,  for a partition $\lambda$ of $n$, we use $f_\lambda(k)$ to denote the number of appearance of $k$ in $\lambda$. Using this notation, we may write $\lambda$ as $(n^{f_\lambda(n)},(n-1)^{f_\lambda(n-1)},\ldots,1^{f_\lambda(1)})$, where when $f_\lambda(i)=0$, we may omit the term $i^0$, and we write $i^1$ as $i$ for simplicity. For example, $\lambda=(5,3,3,1,1,1)$. We have that $f_\lambda(5)=1$, $f_\lambda(3)=2$, $f_\lambda(1)=3$ and $f_\lambda(i)=0$ for all $i\ne 1,3,5$. We can write $\lambda=(5,3^2,1^3)$ instead of $(5,3,3,1,1,1)$.
Let $\lambda,\mu$ be two partitions. We use $\lambda\cup\mu$ to denote the partition
\[(m^{f_\lambda(m)+f_\mu(m)},(m-1)^{f_\lambda(m-1)+f_\mu(m-1)},\ldots,1^{f_\lambda(1)+f_\mu(1)}),\]where $m=\max\{\lambda_1,\mu_1\}$.
Moreover, if $f_\lambda(k)\ge f_\mu(k)$ for all $k$, then we denote the partition
\[(n^{f_\lambda(n)-f_\mu(n)},(n-1)^{f_\lambda(n-1)-f_\mu(n-1)},\ldots,1^{f_\lambda(1)-f_\mu(1)})\]
by $\lambda\setminus\mu$. For example, $\lambda=(6,5^2,2^4,1^5)$ and $\mu=(5,2^3,1^2)$, we have $\lambda\cup\mu=(6,5^3,2^7,1^7)$ and $\lambda\setminus\mu=(6,5,2,1^3)$.

To give the generating function of $b_{t,2}(n)-b_{t,1}(n)$, we need to recall the generating functions of $b_{t,1}(n)$ and $b_{t,2}(n)$, which were given by Singh and Barman \cite{Singh-Barman-2024} and Kim \cite{Kim} as follows.
\begin{thm}[\cite{Singh-Barman-2024}]
  For $t \geq 2$, we have
  \begin{align}\label{equ-bt1}
    \sum_{n=0}^\infty b_{t,1}(n)q^n=\frac{(q^t;q^t)_\infty}{(q;q)_\infty}\left(\frac{q}{1-q}-\frac{q^t}{1-q^t}\right).
  \end{align}
\end{thm}

\begin{thm}[\cite{Kim}]
    For $t\geq 2$, we have
    \begin{equation}\label{equ-bt2}
\sum_{n=0}^\infty b_{t,2}(n)q^n =\frac{(q^t;q^t)_\infty}{(q;q)_\infty}\left(\frac{2q^2}{1-q^2}-\frac{q^t}{1-q^t}+\frac{q^{2t-1}-q^{2t}+q^{2t+1}}{1-q^{2t}}\right).
    \end{equation}
\end{thm}

Here we use the standard $q$-series notation
\[(a;q)_n=\prod_{i=0}^{n-1}(1-aq^i)\quad\text{and}\quad (a;q)_\infty=\prod_{i=0}^{\infty}(1-aq^i)\]

We first establish the following generating function on the difference between $b_{t,2}(n) $ and $ b_{t,1}(n) $.

\begin{thm}\label{-A+B+C}
    For $t \geq 2 $, we have
    \begin{equation}
        \sum_{n=0}^{\infty}\left(b_{t,2}(n)-b_{t,1}(n)\right)q^{n}=(-A+B+C).
    \end{equation}
    where $A=\dfrac{q}{1-q^2} \dfrac{(q^t;q^t)_\infty}{(q^2;q)_\infty}$, $B=\dfrac{q^{2t-1}}{1-q^{2t}}\dfrac{(q^t;q^t)_\infty}{(q^2;q)_\infty}$, and $C=\dfrac{q^{2t+1}}{1-q^{2t}}\dfrac{(q^t;q^t)_\infty}{(q;q)_\infty}$
\end{thm}
\pf From \eqref{equ-bt1} and \eqref{equ-bt2}, we deduce that
\begin{equation}
\begin{split}
 &\sum_{n=0}^{\infty} (b_{t,2}(n)-b_{t,1}(n))q^n\\
 =&\frac{(q^t;q^t)_\infty}{(q;q)_\infty}\left(\frac{2q^2}{1-q^2}-\frac{q^t}{1-q^t}+\frac{q^{2t-1}-q^{2t}+q^{2t+1}}{1-q^{2t}}-(\frac{q}{1-q}-\frac{q^t}{1-q^t})\right)\\
  =&\frac{(q^t;q^t)_\infty}{(q;q)_\infty}\left(\frac{2q^2-q(1+q)}{1-q^2}+\frac{q^{2t-1}(1-q)}{1-q^{2t}}+\frac{q^{2t+1}}{1-q^{2t}}\right)\\
  =&\frac{(q^t;q^t)_\infty}{(q;q)_\infty}\frac{-q(1-q)}{1-q^2}+\frac{(q^t;q^t)_\infty}{(q;q)_\infty}\frac{q^{2t-1}(1-q)}{1-q^{2t}}+\frac{(q^t;q^t)_\infty}{(q;q)_\infty}\frac{q^{2t+1}}{1-q^{2t}}\\
  =&-\frac{(q^t;q^t)_\infty}{(q^2;q)_\infty}\frac{q}{1-q^2}+\frac{(q^t;q^t)_\infty}{(q^2;q)_\infty}\frac{q^{2t-1}}{1-q^{2t}}+\frac{(q^t;q^t)_\infty}{(q;q)_\infty}\frac{q^{2t+1}}{1-q^{2t}}\\
  =&(-A+B+C).
\end{split}
\end{equation}
\qed

Clearly, $A$, $B$ and $C$ have non-negative coefficients. We proceed to show that $-A+C$ has no negative coefficients for $n\geq 192t^5-192t^4-24t^3+24t^2+6t+2$.

\begin{thm}\label{-A+C}
    Given $t\ge 2$ , we have that the coefficient of $q^n$ in
    \begin{equation}
        (-A+C)=-\frac{q}{1-q^2} \frac{(q^t;q^t)_\infty}{(q^2;q)_\infty}+\frac{q^{2t+1}}{1-q^{2t}}\frac{(q^t;q^t)_\infty}{(q;q)_\infty}.
    \end{equation}
    are non-negative for  $n\geq 192t^5-192t^4-24t^3+24t^2+6t+2$.
\end{thm}

In order to show Theorem \ref{-A+C}, we will give combinatorial interpretations to $A$ and $C$ respectively.    Let $O_t(n)$ denote the set of $t$-regular partitions of $n$ in which $1$ appears an odd number of times. Clearly
\[A=\sum_{n=0}^\infty \# O_t(n) q^n.\]
Let $R_t(n)$ denote the set of the partition of $n$ which $kt$ can not appear as a part except for $k=2$. Moreover, the part $2t+1$ must appear at least once. Clearly,
\[C=\sum_{n=0}^\infty \# R_t(n) q^n.\]

Thus Theorem \ref{-A+C} is equivalent to the following Theorem.

\begin{thm}\label{thm-123}
    For $n\ge 192t^5-192t^4-24t^3+24t^2+6t+2$, there is an injection $\Phi$ from $O_t(n)$ into $R_t(n)$.
\end{thm}

To prove Theorem \ref{thm-123}, we first partition the set $O_t(n)$ into five disjoint subsets as follows.
\begin{itemize}
    \item [(1)] $O_t^1(n)$ is the set of partitions of $\lambda \in O_t(n)$ such that there exists $k\ge 1$ such that $2kt+1$ is a part.
    \item [(2)] $O_t^2(n)$ is the set of partitions of $\lambda \in O_t(n)$ such that $\lambda_i \not\equiv 1 \pmod{2t} $ for any $1\le i\le \ell(\lambda)$ and $\lambda_i>1$. Moreover $\lambda_1 \geq 8t^2+1$.
    \item [(3)] $O_t^3(n)$ is the set of partitions of $\lambda \in O_t(n)$ such that $\lambda_i \not\equiv 1 \pmod{2t} $ for any $1\le i\le \ell(\lambda)$ and $\lambda_i>1$. Moreover $\lambda_1 < 8t^2+1$. And there exists $l \geq 2$ as a part that $f_\lambda(l) \geq 6t+1$.
    \item[(4)] $O_t^4(n)$ is the set of partitions of $\lambda \in O_t(n)$ such that $\lambda_i \not\equiv 1 \pmod{2t} $ for any $1\le i\le \ell(\lambda)$ and $\lambda_i>1$. Moreover $\lambda_1 < 8t^2+1$. For any $\lambda_i \geq 2$, $f_\lambda(\lambda_i)<6t+1$. And $f_\lambda(1)\geq 12t+3$.
    \item[(5)] $O_t^5(n)$ is the set of partitions of $\lambda \in O_t(n)$ such that $\lambda_i \not\equiv 1 \pmod{2t} $ for any $1\le i\le \ell(\lambda)$ and $\lambda_i>1$. Moreover $\lambda_1 < 8t^2+1$. For any $\lambda_i \geq 2$, $f_\lambda(\lambda_i)<6t+1$. And $f_\lambda(1)\leq 12t+2$.
\end{itemize}

We next list $4$ disjoint subsets of $R_t(n)$ as given below.
\begin{itemize}
    \item [(1)] $R_t^1(n)$ is the set of partitions $\mu \in R_t(n)$ such that $f_\mu(1) \equiv 1 \pmod{2}$.
    \item [(2)] $R_t^2(n)$ is the set of partitions $\mu \in R_t(n)$ such that $f_\mu(1) \equiv 0 \pmod{2}$. Moreover, $f_\mu(4t+1)+f_\mu(2t+1)\ge 2$ and for $k\geq3$, $f_\mu(2kt+1)=0$.
    \item [(3)] $R_t^3(n)$ is the set of partitions $\mu \in R_t(n)$ such that $f_\mu(1) \equiv 0 \pmod{2}$. Moreover, $f_\mu(6t+1)>0$ and $f_\mu(2kt+1)=0$ for $k\neq 1$ or $3$.
    \item [(4)] $R_t^4(n)$ is the set of partitions $\lambda \in R_t(n)$ such that $f_\mu(1) \equiv 0 \pmod{2}$. Moreover, $f_\mu(8t+1)>0$ and $f_\mu(2kt+1)=0$ if $k\not= 1$ or $4$.
    \end{itemize}

To describe the injection $\Phi$, we proceed to build four injections $\phi_i$ from $O_t(n)$ into $R_t^i(n)$ for $1\le i\le 4$. The injection $\phi_1$ is described as follows.

\begin{lem}\label{lem-phi1}
    There is a bijection $\phi_1$ from $O_t^1(n)$ into $R_t^1(n)$.
\end{lem}

\pf Given $\lambda \in O_t^1(n)$, by definition, we see that there exists $k\ge 1$ such that $2kt+1$ is a part. We choose such $k$ to be minimum. Define
\begin{equation}
    \mu=\phi_1(\lambda)=\lambda \setminus (2kt+1) \cup (2t+1,2t^{k-1})
\end{equation}

We proceed to verify that $\mu \in R_t^1(n)$. Since $\lambda$ is a $t$-regular partition in which $1$ appears an odd number of times, we see that the part $2t+1$ is a part of the partition $\mu$ and $f_\mu(kt)= 0$ for any $k\ne 2$. Moreover, $f_\mu(1) =f_{\lambda}\equiv 1 \pmod{2}$. Furthermore,
\begin{equation}
    |\mu|=|\lambda|-|2kt+1|+|2t+1|+(k-1)|2t|=n.
\end{equation}
This yields $\mu \in R_t^1(n)$.

We next show that $\phi_1$ is a bijection.   We proceed to establish the inverse map $\phi_1^{-1}$. For any $\mu \in R_t^1(n)$, by definition, we see that $f_\mu(1) \equiv 1 \pmod{2}$. Moreover,   $2t+1$ appears in $\mu$ and $f_\mu(2t)\geq 0$. Set $k=1+{f_\mu(2t)}$. Define
\begin{equation}
    \phi^{-1}_1(\mu)=\mu\setminus (2t+1,2t^{f_\mu(2t)})\cup(2kt+1).
\end{equation}
It is routine to check that $\phi_1^{-1}$ is the inverse map of $\phi_1$. Thus $\phi_1$ is a bijection.\qed

For example, set $t=4$, $n=65$, $\lambda=(17,15,13,10,7,5,3,2,1,1,1)$. Clearly $\lambda \in O_4^1(65)$. Applying  $\phi_1$ on $\lambda$, we see that  $k=2$ and
\begin{equation}
    \mu=\phi_1(\lambda)=(15,13,10,9,8,5,3,2,1,1,1).
\end{equation}
So $\mu \in R_4^1(n)$. Applying $\phi_1^{-1}$ on $\mu$, we recover $\lambda$.

We next give the injection $\phi_2$.

\begin{lem}
    There is an injection $\phi_2$ from $O_t^2(n)$ into $R_t^2(n)$.
\end{lem}
\pf Given $\lambda \in O_t^2(n)$, by definition, we see that $\lambda_i \not\equiv 1 \pmod{2t} $ for any $\lambda_i\geq 2$ and $1\le i\le \ell(\lambda)$. Moreover, $\lambda_1 \geq 8t^2+1$. Notice that $(2t+1,4t+1)=1$, from $\lambda_1\ge 8t^2+1$, we see that there exists unique integer $x,y$ such that $\lambda_1-1=x(2t+1)+y(4t+1)$, where $x\geq 0$ and $0\leq y\leq2t$. There are two cases.

Case 1. If $x\ne 0$, then define
\begin{equation}
    \mu= \phi_2(n)= \lambda\setminus (\lambda_1)\cup ((4t+1)^y,(2t+1)^x,1).
\end{equation}

Case 2. If $x=0$,  then define
\begin{equation}
    \mu= \phi_2(n)= \lambda\setminus (\lambda_1)\cup ((4t+1)^{y-1},2t+1,2t,1).
\end{equation}

We proceed to verify that $\mu \in R_t^2(n)$. Since $\lambda$ is a $t$-regular partition in which $1$ appears an odd number of times, we see that $f_\mu(1)=f_\lambda(1)+1\equiv 0 \pmod{2}$ and $2t+1$ is a part of $\mu$. Moreover, note that $x+y\ge \frac{(8t^2+1)}{4t+1}\ge 2t-1$, which implies $f_\mu(2t+1)+f_\mu(4t+1)=x+y\ge 2$ in either case, and  $f_\mu(2kt+1)=0$ for $k\geq3$. Furthermore,  when $x\ne 0$,
\begin{equation}
    |\mu|=|\lambda|-|\lambda_1|+x|2t+1|+y|4t+1|+1=n.
\end{equation}
And when $x=0$,
\begin{equation}
    |\mu|=|\lambda|-|\lambda_1|+(y-1)|4t+1|+|2t+1|+|2t|+1=n.
\end{equation}
This yields $\mu \in R_t^2(n)$.

We next show that $\phi_2$ is an injection. To this end, let $I_2(n)$ denote the image set of $\phi_2$, which has been shown to be a subset of $R_t^2(n)$. We proceed to establish a map $\psi_2$  from $I_2(n)$ to $O_t^2(n)$, such that for any $\lambda \in O_t^2(n)$,
\begin{equation}\label{equ-psi-2}
    \psi_2(\phi_2(\lambda))=(\lambda).
\end{equation}
This implies $\phi_2$ is an injection.

We now describe $\psi_2$. For any $\mu \in R_t^2(n)$, by the definition of $R_t^2(n)$,  we see that $f_\mu(1) \equiv 0 \pmod{2}$,  $f_\mu(2t+1)\ge1$  and $f_\mu(2t+1)+f_\mu(4t+1)\ge2$. Moreover, from the construction of $\phi_2$, we see that  $f_\mu(2t)=0$ or $1$, $f_\mu(1)\ge 1$. Let
\[l=1+2tf_\mu(2t)+(2t+1)f_\mu(2t+1)+(4t+1)f_\mu(4t+1),\]
then we have both $l\ge \mu_2$ and $l\ge 8t^2+1$.
Define
\begin{equation}
    \lambda=\psi_2(\mu)=\mu \setminus ((4t+1)^{f_\mu(4t+1)},(2t+1)^{f_\mu(2t+1)},(2t)^{f_\mu(2t)},1)\cup (l).
\end{equation}
It is routine to check that $\psi_2(\mu) \in O_t^2(n)$ and \eqref{equ-psi-2} holds. Thus $\phi_2$ is an injection.\qed

For example, set $t=4$, $n=221$ and $\lambda=(137,33,29,11,5,3,1,1,1)$. Clearly $\lambda \in O_4^2(221)$. Applying $\phi_2$ on $\lambda$, we see that $x=0$, $y=7$. Hence
\begin{equation}
    \mu=\phi_2(\lambda)=(33,29,17^7,11,9,8,5,3,1^4)
\end{equation}
So $\mu \in R_4^2(n)$. Applying $\psi_2$ on $\mu$, we recover $\lambda$.

For another example, set $t=4$, $n=242$ and $\lambda=(157,34,29,11,5,3,1,1,1)$. Clearly $\lambda \in O_4^2(242)$. Applying $\phi_2$ on $\lambda$, we see that $x=6$, $y=6$. Hence
\begin{equation}
    \mu=\phi_2(\lambda)=(34,29,17^6,11,9^6,5,3,1^4)
\end{equation}
So $\mu \in R_4^2(n)$. Applying $\psi_2$ on $\mu$, we recover $\lambda$.

We next give the injection $\phi_3$.

\begin{lem}
    There is an injection $\phi_3$ from $O_t^3(n)$ into $R_t^3(n)$.
\end{lem}
\pf Given $\lambda \in O_t^3(n)$, $\lambda_i \not\equiv 1 \pmod{2t} $ for any $\lambda_i\geq 2$ and $1\le i\le \ell(\lambda)$. Moreover $\lambda_1 < 8t^2+1$, and there exists $l \geq 2$ as a part that $f_\lambda(l) \geq 6t+1$. We choose such $l$ to be minimum. Define
\begin{equation}
    \mu=\phi_3(\lambda)=\lambda\setminus(l^{6t+1})\cup((6t+1)^{l-1},(2t+1)^2,1^{2t-1}).
\end{equation}

We proceed to verify that $\mu \in R_t^3(n)$. Since $\lambda$ is a $t$-regular partition in which $1$ appears an odd number of times, we see that $f_\mu(1)=f_\lambda(1)+2t-1\equiv 0 \pmod{2}$ and $2t+1$ is a part of $\mu$. Moreover, $f_\mu(6t+1) \geq 1$  and $f_\mu(2kt+1)=0$ if $k\not= 1$ or $3$. Furthermore,
\begin{equation}
    |\mu|=|\lambda|-(6t+1)l+(l-1)(6t+1)+2(2t+1)+(2t-1)=n
\end{equation}
This yields $\mu \in R_t^3(n)$.

We next show that $\phi_3$ is an injection. To this end, let $I_3(n)$ denote the image set of $\phi_3$, which has been shown to be a subset of $R_t^3(n)$. We proceed to establish a map $\psi_3$  from $I_3(n)$ to $O_t^3(n)$, such that for any $\lambda \in O_t^3(n)$,
\begin{equation}\label{equ-psi-3}
    \psi_3(\phi_3(\lambda))=(\lambda).
\end{equation}
This implies $\phi_3$ is an injection.

We now describe $\psi_3$. For any $\mu \in R_t^3(n)$, by the definition of $R_t^3(n)$, we see that $f_\mu(1) \equiv 0 \pmod{2}$, $f_\mu(6t+1)\geq 1$ and $f_\mu(2kt+1)=0$ except $k=1$ or $3$. Moreover,from the construction of $\phi_3$, we see that $f_\mu(2t+1)=2$ and $f_\mu(1) \geq2t-1$.  Define
\begin{equation}
    \lambda=\psi_3(\mu)=\mu\setminus((6t+1)^{f_\mu(6t+1)},(2t+1)^2,1^{2t-1})\cup((f_\mu(6t+1)+1)^{6t+1})
\end{equation}
It is routine to check that $\psi_3(\mu) \in O_t^3(n)$ and \eqref{equ-psi-3} holds. Thus $\phi_3$ is an injection.\qed

For example, set $t=4$, $n=128$, $n=(17,13,11,9,3^{25},1,1,1)$. Clearly $\lambda \in O_4^3(128)$.
\begin{equation}
    \mu=\phi_3(\lambda)=(25,25,17,13,11,9,9,9,1^{10})
\end{equation}
So $\mu \in R_4^3(n)$. Applying $\psi_3$ on $\mu$, we recover $\lambda$.

We next give the injection $\phi_4$.

\begin{lem}\label{lem-phi4}
    There is an injection $\phi_4$ from $O_t^4(n)$ into $R_t^4(n)$.
\end{lem}

\pf Given $\lambda \in O_t^4(n)$, $\lambda_i \not\equiv 1 \pmod{2t} $ and $f_\lambda(\lambda_i)\leq 6t+1$ for any $\lambda_i \ge2$. Moreover $\lambda_1 < 8t^2+1$ and $f_\lambda(1)\geq 12t+3$. Define
\begin{equation}
    \mu=\phi_4(\lambda)=\lambda\setminus(1^{12t+3})\cup(8t+1,(2t+1)^2).
\end{equation}

We proceed to verify that $\mu \in R_t^4(n)$. Since $\lambda$ is a $t$-regular partition in which $1$ appears an odd number of times, we see that $f_\mu(1)=f_\lambda(1)-(12t+3)\equiv 0 \pmod{2}$ and $2t+1$ is a part of $\mu$.
Moreover, $f_\mu(8t+1) = 1$, $f_\mu(2t+1)=2$ and $f_\mu(2kt+1)=0$ if $k\not= 1$ or $4$. Furthermore,
\begin{equation}
    |\mu|=|\lambda|-(12t+3)+(8t+1)+2(2t+1)=n.
\end{equation}
This yields $\mu \in R_t^4(n)$.

We next show that $\phi_4$ is an injection. To this end, let $I_4(n)$ denote the image set of $\phi_4$, which has been shown to be a subset of $R_t^4(n)$. We proceed to establish a map $\psi_4$  from $I_4(n)$ to $O_t^4(n)$, such that for any $\lambda \in O_t^4(n)$,
\begin{equation}\label{equ-psi-4}
    \psi_4(\phi_4(\lambda))=(\lambda).
\end{equation}
This implies $\phi_4$ is an injection.

We now describe $\psi_4$. For any $\mu \in R_t^4(n)$, by the definition of $R_t^4(n)$, we see that $f_\mu(1) \equiv 0 \pmod{2}$. Moreover, from the construction of $\phi_4$,  $f_\mu(2t+1)=2$, $f_\mu(8t+1)= 1$ and $f_\mu(2kt+1)=0$ except $k=1$ or $4$. Define
\begin{equation}
    \lambda=\psi_4(\mu)=\mu\setminus(8t+1,(2t+1)^2)\cup(1^{12t+3})
\end{equation}
It is routine to check that $\psi_4(\mu) \in O_t^4(n)$ and \eqref{equ-psi-4} holds. Thus $\phi_4$ is an injection.\qed

For example, set $t=4$, $n=85$, $\lambda=(13,7,6,2,2,1^{55})$. Clearly, $\lambda \in O_4^4(85)$.
\begin{equation}
    \mu=\phi_4(\lambda)=(33,13,9,9,7,6,2,2,1,1,1,1)
\end{equation}

So $\mu \in R_4^4(n)$. Applying $\psi_4$ on $\mu$, we recover $\lambda$.

We next show that when  $n\geq 192t^5-192t^4-24t^3+24t^2+6t+2$, $O_5^t(n)$ is empty. In fact, if $\lambda \in O_t^5(n)$, then
\begin{equation}
    |\lambda| \leq 6t\sum_{n=2}^{8t^2}n -6t\sum_{n=1}^{8t}nt+(12t+1)=192t^5-192t^4-24t^3+24t^2+6t+1,
\end{equation}
 which is contradict to  $n\geq 192t^5-192t^4-24t^3+24t^2+6t+2$.

 We are now in a position to prove Theorem \ref{thm-123}.

{\noindent \it Proof of Theorem \ref{thm-123}.} For any given integer $t$, if $n\geq 192t^5-192t^4-24t^3+24t^2+6t+2$, $\lambda \in O_t(n)$, define
\begin{equation}
    \Phi(\lambda)=\begin{cases}
    \phi_1(\lambda),&\text{if }\lambda \in O_t^1(n);\\
    \phi_2(\lambda),&\text{if }\lambda \in O_t^2(n);\\
    \phi_3(\lambda),&\text{if }\lambda \in O_t^3(n);\\
    \phi_4(\lambda),&\text{if }\lambda \in O_t^4(n).\\
    \end{cases}
\end{equation}
From Lemma \ref{lem-phi1} $\sim$ Lemma \ref{lem-phi4}, we see that $\Phi$ is indeed an injection. \qed

\section{A proof of Theorem \ref{thm-bt2-bt3}}\label{3}

In this section, we give a proof of Theorem \ref{thm-bt2-bt3}.  We first write the  generating function of $b_{t,2}(n)- b_{t,3}(n) $ as a summation of three summands, as stated in Theorem \ref{D+E+F}. We then prove that each summand is non-negative for specific $n$, as stated in Theorem \ref{thm-D}, Theorem \ref{thm-E} and Theorem \ref{thm-F}. In this way, the proof of Theorem \ref{thm-bt2-bt3} are given.

We first recall the generating function of $b_{t,3}(n)$ as given by Kim \cite{Kim}.

\begin{thm}[\cite{Kim}]\label{equ-bt3}
    For an integer $t\geq2$, we have
    \begin{equation}
        \begin{aligned}
\sum_{n=0}^{\infty}b_{t,3}(n)q^n=\frac{(q^t;q^t)_\infty}{(q;q)_\infty}\left(\frac{3q^{3}}{1-q^{3}}-\frac{q^{t}}{1-q^{t}}+\frac{q^{2t-2}-q^{2t}+q^{2t+2}}{1-q^{2t}} \right.\\
\left.-\frac{q^{3t-3}-q^{3t-2}-q^{3t-1}+2q^{3t}-q^{3t+1}-q^{3t+2}+q^{3t+3}}{1-q^{3t}}\right).
    \end{aligned}
    \end{equation}
\end{thm}

We first establish the following generating function on the difference between $b_{t,2}(n) $ and $ b_{t,3}(n) $.

\begin{thm}\label{D+E+F}
For $n\geq2$, we have
    \begin{equation}
    \sum_{n=0}^\infty(b_{t,2}(n)-b_{t,3}(n))q^n=(D+E+F).
\end{equation}
where
$$D=\frac{(q^t;q^t)_\infty}{(q^2;q)_\infty}\frac{q^2+q^4}{1-q^6}-\frac{(q^t;q^t)_\infty}{(q;q)_\infty}\frac{q^{2t-2}(1-q)(1-q^3)}{1-q^{2t}}, $$ $$E=\frac{(q^t;q^t)_\infty}{(q^2;q)_\infty}\frac{q^2-q^3+q^5}{1-q^6}$$
and
$$F=\frac{(q^t;q^t)_\infty}{(q;q)_\infty}\frac{q^{3t-3}(1+q^3)(1-q)(1-q^2)}{1-q^{3t}}.$$
\end{thm}

\pf From \eqref{equ-bt2} and  \eqref{equ-bt3}, we deduce that
\begin{equation*}
    \begin{split}
 &\sum_{n=0}^{\infty} \left(b_{t,2}(n)-b_{t,3}(n)\right)q^n\\
 =&\frac{(q^t;q^t)_\infty}{(q;q)_\infty}\left(\left(\frac{2q^2}{1-q^2}-\frac{q^t}{1-q^t}+\frac{q^{2t-1}-q^{2t}+q^{2t+1}}{1-q^{2t}}\right)\right.\\
 &-\left(\frac{3q^{3}}{1-q^{3}}-\frac{q^{t}}{1-q^{t}}+\frac{q^{2t-2}-q^{2t}+q^{2t+2}}{1-q^{2t}} -\right.\\
 &\left.\left.\frac{q^{3t-3}-q^{3t-2}-q^{3t-1}+2q^{3t}-q^{3t+1}-q^{3t+2}+q^{3t+3}}{1-q^{3t}}\right)\right)\\
 =&\frac{(q^t;q^t)_\infty}{(q;q)_\infty}\left(\frac{2q^2}{1-q^2}-\frac{3q^{3}}{1-q^{3}}-\frac{q^{2t-2}(1-q)(1-q^3)}{1-q^{2t}}+\frac{q^{3t-3}(1+q^3)(1-q)(1-q^2)}{1-q^{3t}}\right)\\
 =&\frac{(q^t;q^t)_\infty}{(q;q)_\infty}\left(\frac{q^2(1-q)^2+(q^2+q^4+q^5)(1-q)}{1-q^6}-\frac{q^{2t-2}(1-q)(1-q^3)}{1-q^{2t}}+\right.\\
 &\left.\frac{q^{3t-3}(1+q^3)(1-q)(1-q^2)}{1-q^{3t}}\right)\\
 =&\frac{(q^t;q^t)_\infty}{(q^2;q)_\infty}\frac{q^2+q^4}{1-q^6}-\frac{(q^t;q^t)_\infty}{(q;q)_\infty}\frac{q^{2t-2}(1-q)(1-q^3)}{1-q^{2t}}+\frac{(q^t;q^t)_\infty}{(q^2;q)_\infty}\frac{q^2-q^3+q^5}{1-q^6}+\\
 &\frac{(q^t;q^t)_\infty}{(q;q)_\infty}\frac{q^{3t-3}(1+q^3)(1-q)(1-q^2)}{1-q^{3t}}.\\
 =&D+E+F
\end{split}
\end{equation*}
\qed
We proceed to study the non-negativity of  $D$, $E$ and $F$. We first give the following non-negativity of $D$ as stated in Theorem \ref{thm-D}.

\begin{thm}\label{thm-D}
    For any $t\ge 2$ and $n\ge 0$, the coefficient of $q^n$ in $D$ is non-negative except for $(t,n)=(2,6)$.
\end{thm}

We first consider the case $t\ge 4$.
In order to show Theorem \ref{thm-D},  we  give combinatorial interpretations to $\frac{(q^t;q^t)_\infty}{(q^2;q)_\infty}\frac{q^2+q^4}{1-q^6}$ and $\frac{(q^t;q^t)_\infty}{(q;q)_\infty}\frac{q^{2t-2}(1-q)(1-q^3)}{1-q^{2t}}$. Let $S_t(n)$ denote the set of $t$-regular partition $\lambda$ of $n$ in which $f_\lambda(1) \equiv2 \pmod{6}$ or $f_\lambda(1) \equiv4 \pmod{6}$. Clearly
\begin{equation}
    \frac{(q^t;q^t)_\infty}{(q^2;q)_\infty}\frac{q^2+q^4}{1-q^6}=\sum_{n=0}^\infty \# S_t(n) q^n.
\end{equation}
Let $A_t(n)$ denote the set of the partition $\lambda$ of $n$ which $f_\lambda(3)=0$ and $f_\lambda(1)\equiv-2 \pmod{2t}$. Clearly,
\begin{equation}
    \frac{(q^t;q^t)_\infty}{(q;q)_\infty}\frac{q^{2t-2}(1-q)(1-q^3)}{1-q^{2t}}=\sum_{n=0}^\infty \# A_t(n) q^n.
\end{equation}

Thus when $t\ge 4,$ Theorem \ref{thm-D} is equivalent to the following theorem.

\begin{thm}\label{thm-A-S}
    For  $t\geq4$, there is an injection $\Gamma$ from $A_t(n)$ into $S_t(n)$.
\end{thm}

To prove Theorem \ref{thm-A-S}, we first partition the set $A_t(n)$ into three disjoint subsets as follows.

\begin{itemize}
    \item[(1)] $A_t^1(n)$ is the set of partitions of $\lambda \in A_t(n)$ such that $f_\lambda(1) \equiv 2 \pmod{6}$.
    \item[(2)] $A_t^2(n)$ is the set of partitions of $\lambda \in A_t(n)$ such that $f_\lambda(1) \equiv 4 \pmod{6}$.
    \item[(3)] $A_t^3(n)$ is the set of partitions of $\lambda \in A_t(n)$ such that $f_\lambda(1) \equiv 0 \pmod{6}$.
\end{itemize}

We next list $3$ disjoint subsets of $S_t(n)$ as given below.

\begin{itemize}
    \item[(1)] $S_t^1(n)$ is the set of partitions $\mu \in S_t(n)$ such that $f_\mu(1) \equiv 2 \pmod{6}$.
    \item[(2)] $S_t^2(n)$ is the set of partitions $\mu \in S_t(n)$ such that $f_\mu(1) \equiv 4 \pmod{6}$ and $f_\mu(1) \equiv-2\pmod{2t}$.
    \item[(3)] $S_t^3(n)$ is the set of partitions $\mu \in S_t(n)$ such that $f_\mu(1) \equiv 4 \pmod{6}$ and $f_\mu(1) \equiv-4\pmod{2t}$.
\end{itemize}

We next build three injections $\gamma_i$ from $A_t^i(n)$ into $S_t^i(n)$, where $1\le i\le 3.$
It is clear that $A_t^1(n)\subseteq S_t^1(n)$ and $A_t^2(n)\subseteq S_t^2(n)$. Hence we may set the map $\gamma_1$ and $\gamma_2$ to be the identity map.
We next describe the infection $\gamma_3$ from $A_t^3(n)$ into $S_t^3(n)$ in the following lemma.

\begin{lem}\label{lem-gamma-3}
    There is an injection $\gamma_3$ from $A_t^3(n)$ into $S_t^3(n)$.
\end{lem}

\pf Given $\lambda \in A_t^3(n)$, by definition $f_\lambda(1) \equiv 0 \pmod{6}$ and $f_\lambda(1)\equiv-2 \pmod{2t}$. Define
\begin{equation}
    \mu=\gamma_3(\lambda)=\lambda\setminus(1^2)\cup(2).
\end{equation}

We proceed to verify that $\lambda \in S_t^3(n)$. Since $\lambda$ is a $t$-regular partition in which $f_\lambda(1) \equiv-2 \pmod{2t}$ and $f_\lambda(1) \equiv 0 \pmod{6}$, we can see that $f_\mu(1)=f_\lambda(1)-2 \equiv-4 \pmod{2t}$ and $f_\mu(1)=f_\lambda(1)-2 \equiv -2 \pmod{6}$. Furthermore,
\begin{equation}
    |\mu|=|\lambda|-2+2=n
\end{equation}
This yields $\mu \in S_t^3(n)$.

We next show that $\gamma_3$ is an injection. To this end, let $I_3^A(n)$ denote the image set of $\gamma_3$, which has been shown to be a subset of $S_t^3(n)$. We proceed to establish a map $\delta_3$  from $I_3^A(n)$ to $A_t^3(n)$, such that for any $\lambda \in A_t^3(n)$,
\begin{equation}\label{equ-delta-3}
    \delta_3(\gamma_3(\lambda))=\lambda.
\end{equation}
This implies $\gamma_3$ is an injection.

We now describe $\delta_3$. For any $\mu \in S_t^3(n)$, by the definition of $S_t^3(n)$, we see that $f_\mu(1)\equiv-4 \pmod{2t}$ and $f_\mu(1)\equiv -2 \pmod{6}$. Moreover, from the construction of $\gamma_3$, we see that $f_\mu(2) \geq 1$ and $f_\mu(3)=0$. Define
\begin{equation}
    \lambda =\delta_3(\mu)=\mu\setminus(2) \cup(1^2).
\end{equation}
It is routine to check that $\gamma_3(\mu)\in A_t^3(n)$ and \ref{equ-delta-3} holds. Thus $\gamma_3$ is an injection. \qed

We are now in a position to prove Theorem \ref{thm-A-S}.

{\noindent \it Proof of Theorem \ref{thm-A-S}.} For any given integer $t\geq4$, define
\begin{equation}
    \Gamma(\lambda)=\begin{cases}
    \gamma_1(\lambda),&\text{if }\lambda \in A_t^1(n);\\
    \gamma_2(\lambda),&\text{if }\lambda \in A_t^2(n);\\
    \gamma_3(\lambda),&\text{if }\lambda \in A_t^3(n).\\
    \end{cases}
\end{equation}
From Lemma \ref{lem-gamma-3}, we see that $\Gamma$ is an injection.\qed

We now give a proof of Theorem \ref{thm-D}.

{\noindent \it Proof of Theorem \ref{thm-D}.} If $t\geq4$, by Theorem \ref{thm-A-S}, $D$ have non-negative coefficients.

Now we proceed to study the positivity of $D$ for $t\le 3$.

When $t=3$, we see that
\begin{equation}
\begin{aligned}
    D=&\frac{(q^3;q^3)_\infty}{(q;q)_\infty}\left(\frac{(q^2+q^4)(1-q)}{1-q^6}-\frac{q^4(1-q)(1-q^3)}{1-q^6}\right)\\
    =&\frac{(q^3;q^3)_\infty}{(q;q)_\infty}\left(\frac{q^2(1-q)}{1-q^6}+\frac{q^7(1-q)}{1-q^6}\right)\\
    =&\frac{(q^3;q^3)_\infty}{(q^2;q)_\infty}\frac{q^2+q^7}{1-q^6}.
\end{aligned}
\end{equation}
It is clear that $D$ has non-negative coefficients when $t=3$.

When $t=2$, we aim to show that $D$ has non-negative coefficients of $q^n$ except for $n=6$. In this case,
    \begin{align}\label{equ-D123}
        D=&\frac{(q^2;q^2)_\infty}{(q;q)_\infty}\left(\frac{(q^2+q^4)(1-q)}{1-q^6}-\frac{q^2(1-q)(1-q^3)}{1-q^4}\right)\nonumber\\
        =&\frac{(q^2;q^2)_\infty}{(q^2;q)_\infty}\frac{(q^2+q^4)(1+q^6)-q^2(1-q^3)(1+q^4+q^8)}{1-q^{12}}\nonumber\\
        =&\frac{(q^2;q^2)_\infty}{(q^2;q)_\infty}\frac{q^4-q^6}{1-q^{12}}+\frac{(q^2;q^2)_\infty}{(q^2;q)_\infty}\frac{q^5+q^8+q^9+q^{13}}{1-q^{12}}
    \end{align}
It is clear that the second summation in \eqref{equ-D123} has non-negative coefficients. We proceed to prove that
\[\frac{(q^2;q^2)_\infty}{(q^2;q)_\infty}\frac{q^4-q^6}{1-q^{12}}\]
also has non-negative coefficients of $q^n$ except for $n=6$.  Let $D_1(n)$ denote the set of $2$-regular partition $\lambda$ of $n$ in which $f_\lambda(1)\equiv4\pmod{12} $. Clearly
\begin{equation}
    \sum_{n=0}^\infty\#D_1(n)q^n=\frac{(q^2;q^2)_\infty}{(q^2;q)_\infty}\frac{q^4}{1-q^{12}}.
\end{equation}
Let $D_2(n)$ denote the set of $2$-regular partition $\lambda$ of $n$ in which $f_\lambda(1)\equiv6\pmod{12} $. Clearly
\begin{equation}
\sum_{n=0}^\infty\#D_2(n)q^n=\frac{(q^2;q^2)_\infty}{(q^2;q)_\infty}\frac{q^6}{1-q^{12}}.
\end{equation}

Thus, to show Theorem \ref{thm-D}, it suffices to build an injection $\epsilon$ from $D_2(n)$ into $D_1(n)$ for any integer $n\geq7$.
Given $\lambda\in D_2(n)$, there are two cases:

Case 1: $\lambda_1\ge 3$. In this case, define
\begin{equation}
    \mu=\epsilon(\lambda)=\lambda\setminus(\lambda_1,1^2)\cup(\lambda_1+2)
\end{equation}
We proceed to verify that $\mu \in D_1(n)$. Since $\lambda$ is $2$-regular partition in which $f_\lambda(1)\equiv6 \pmod{12}$, we see that $f_\mu(1)=f_\lambda(1)-2\equiv 4\pmod{12}$ and $\lambda_1+2\equiv1\pmod{2}$. Moreover,
\begin{equation}
    |\mu|=|\lambda|-\lambda_1-2+\lambda_1+2=n
\end{equation}
This yields $\mu \in D_1(n)$.

Case 2: $\lambda_1=1$. In this case, there is only one partition, namely $(1^{12k+6})$ for some $k\ge 1$. Define
\begin{equation}
    \mu=\epsilon(\lambda)=((6k+1)^2,1^4).
\end{equation}

Clearly in Case 1, we have $\mu_1>\mu_2$, and in Case 2, we have $\mu_1=\mu_2$. Thus the image set of these two cases are non-intersect. Moreover,
It is easy to check $\epsilon$ is an injection.

From the above analysis, we prove Theorem \ref{thm-D}.\qed

Now we give the positivity of $E$ as stated in Theorem \ref{thm-E}.

\begin{thm}\label{thm-E}
    For $n\geq 4$ except for $t=2$ and $n=9$, the coefficient of $q^n$ in $E$ is non-negative.
\end{thm}

\pf Clearly,
\begin{equation}
    E=\frac{(q^t;q^t)_\infty}{(q^2;q)_\infty}\frac{q^2-q^3+q^5}{1-q^6}=E_1-E_2
\end{equation}
where $E_1=\frac{(q^t;q^t)_\infty}{(q^2;q)_\infty}\frac{q^2+q^5}{1-q^6}$ and $E_2=\frac{(q^t;q^t)_\infty}{(q^2;q)_\infty}\frac{q^3}{1-q^6}$.

We first consider the case for $t\geq3$.
In order to show Theorem \ref{thm-E}, we give combinatorial interpretations to $E_1$ and $E_2$ respectively. Let $B_t(n)$ denote the set of $t$-regular partition $\lambda$ of $n$ in which either $f_\lambda(1)\equiv2 \pmod{6}$ or $f_\lambda(1) \equiv 5 \pmod{6}$. Clearly
\begin{equation}
    E_1=\sum_{n=0}^\infty\#B_t(n)q^n.
\end{equation}
Let $C_t(n)$ denote the set of $t$-regular partition $\lambda$ of $n$ in which $f_\lambda(1) \equiv 3 \pmod{6}$. Clearly,
\begin{equation}
    E_2=\sum_{n=0}^\infty\#C_t(n)q^n.
\end{equation}

 So it sufficient to build an injection $\tau$ from $C_t(n)$ into $B_t(n)$.

 Given $\lambda \in C_t(n)$, there are three cases.

 Case 1: $2\in \lambda $. Define that
 \begin{equation}
      \mu=\tau(\lambda)=\lambda\setminus(2)\cup(1^2).
 \end{equation}
 It is easy to check $f_\mu(1)\equiv 5 \pmod{6}$. So $\mu\in B_t(n)$.

 Case 2: $2\notin \lambda$, $\lambda_1\geq2$ and $\lambda_1 \not\equiv -1 \pmod{t}$. Define that
 \begin{equation}
     \mu=\tau(\lambda)= \lambda \setminus(\lambda_1,1)\cup(\lambda_1+1)
 \end{equation}
 It is easy to check $f_\mu(1)\equiv 2 \pmod{6}$ and $\lambda_1 \not \equiv 0 \pmod{t}$. So $\mu\in B_t(n)$.

 Case 3: $2\notin \lambda$, $\lambda_1\geq2$ and $\lambda_1 \equiv -1 \pmod{t}$. Define that
 \begin{equation}
     \mu=\tau(\lambda)=\lambda \setminus(\lambda_1,1)\cup(\lambda_1-1,2).
 \end{equation}
 It is easy to check $f_\mu(1)\equiv2 \pmod{6}$ and  $\lambda_1-1 \equiv -2 \pmod{t}$. So $\mu\in B_t(n)$.

 Case 4: $\lambda_1=1$, and $f_\lambda(1)\equiv 3 \pmod{6}$. In this case, $\lambda=(1^{6k+3})$ for some $k\geq 1$. Define
 \begin{equation}
     \mu=\tau(\lambda)=\begin{cases}
         (3,2^{3k-1},1^2),&\text{if } t\geq5;\\
         (5,2^{3k-2},1^2),&\text{if } t=3 \text{ or } 4.
     \end{cases}
 \end{equation}

We first show that the image sets of Case 1, Case 2, Case 3 and Case 4 are pairwise disjoint. By the definition of Case 1, we see that $f_\mu(1)\equiv 5 \pmod{6}$. And $f_\mu(1)\equiv 2\pmod{6}$ in Case 2, Case 3 and Case 4. Moreover, we see that $2\notin \mu$ in Case 2, but $2 \in \mu$ in Case 3 and Case 4. Furthermore, in Case 3, we see that
\begin{equation}
    f_\mu(2)=\begin{cases}
        2&\text{if }(t,\lambda_1)=(4,3);\\
        1&\text{otherwise}.
    \end{cases}
\end{equation}So when $t\geq 5$, $f_\mu(2)=1$ in Case 3 and $f_\mu(2)\geq 2$ in Case 4, which implies the image sets of Case 3 and Case 4 are non-intersect. When $t=3$, since  $4\notin \mu$, we see that $(5,2^{3k-2},1^2)$ cannot be the image of Case 3. When $t=4$,   $(5,2^{3k-2},1^2)$ is not a element of the image of Case 3 because $5 \not\equiv -1$ or $-2\pmod{4}$. From the above analysis, we find that the image sets of Case 1, Case 2, Case 3 and Case 4 are pairwise disjoint.

We next show that in each case, the map $\tau$ is an injection. Let $I'(n)$ denote the image set of $\tau$, which has been shown to be a subset of $B_t(n)$. It sufficient to construct the map $\eta\colon I'(n)\rightarrow C_t(n)$, such that for any $\lambda\in C_t(n)$,
\begin{equation}\label{equ-eta-tau}
    \eta(\tau(\lambda))=\lambda.
\end{equation}

For any $\mu\in I'(n)$, we proceed to illustrate the map $\eta$. There are four cases.

Case 1: $\mu\in B_t(n)$ and $f_\mu(1)\equiv5\pmod{6}$. Define that
\begin{equation}
    \lambda=\eta(\mu)=\mu\setminus(1^2)\cup(2).
\end{equation}

Case 2: $\mu\in B_t(n)$, $f_\mu(1)\equiv2\pmod{6}$, $\mu_1>\mu_2$ and $2\notin\mu$. Define that
\begin{equation}
    \lambda=\eta(\mu)=\mu\setminus(\mu_1)\cup(\mu_1-1,1).
\end{equation}

Case 3: $\mu \in B_t(n)$, $f_\mu(1)\equiv2\pmod{6}$, and $f_\mu(2)=1$. Moreover, $\mu_1 \equiv -1$ or $-2 \pmod{t}$. If $\mu_1 \equiv -1\pmod{t}$, then $f_\mu(\mu_1-1)\geq1$. Let $i$ be the minimum integer such that  $\mu_{i}\equiv-2 \pmod{t}$. Define
\begin{equation}
    \lambda=\eta(\mu)=\mu\setminus(\mu_{i},2)\cup(\mu_{i}+1,1).
\end{equation}

Case 4: $\mu \in B_t(n)$, and
\begin{equation}
    \mu=\begin{cases}
         (3,2^{3k-1},1^2),&\text{if } t\geq5;\\
         (5,2^{3k-2},1^2),&\text{if } t=3 \text{ or } 4.
         \end{cases}
\end{equation}
Define that
\begin{equation}
    \lambda=\eta(\mu)=(1^{|\mu|}).
\end{equation}

It is routine to check that \eqref{equ-eta-tau} holds for any $\lambda\in C_t(n)$. This yields $\tau$ is an injection.

The injection $\tau$ gives a proof of the Theorem \ref{thm-E} for $t\geq3$. Now we prove  Theorem \ref{thm-E} for $t=2$.

From \cite{Weisstein}, we have that for $n\ge 4$,
\begin{equation}
    Q(n)\leq \frac{1}{2}\left(Q(n-1)+Q(n+1)\right)
\end{equation}
where $Q(n)$ denotes the number of the distinct partition of $n$. It is well known that the generating function of $Q(n)$ equals the generating function of the number of $2$-regular partition of $n$.

So that, we have that
\begin{equation}
    \frac{(q^2-q^3)(1-q)(q^2;q^2)_\infty}{(q;q)_\infty}
\end{equation}
has non-negative coefficients for $n\geq 7$. Moreover, by the direct calculation, we have
\begin{equation}
    \frac{(q^2-q^3)(1-q)(q^2;q^2)_\infty}{(q;q)_\infty}=-q^3-q^6+q^{24}+q^{27}+G(q)
\end{equation}
where $G(q)$ has non-negative coefficients. So
\begin{align}\label{equ-G(q)}
    \frac{(q^2-q^3)(1-q)(q^2;q^2)_\infty}{(q;q)_\infty(1-q^6)}
    =&\frac{-q^3-q^6+q^{24}+q^{27}}{1-q^6}+\frac{G(q)}{1-q^6}\nonumber\\
    =&-q^3-q^6-q^9-q^{12}-q^{15}-q^{18}-q^{21}+\frac{G(q)}{1-q^6}
\end{align}
has non-negative coefficients for $n\geq 22$. After checking $4\leq n\leq21$, we have that when $t=2$, $E$ has non-negative coefficients except for $n=3$ or $9$. Combining the case $t\ge 3$, we complete the entire proof of Theorem \ref{thm-E}.\qed

Now we give the positivity of $F$ as stated in Theorem \ref{thm-F}.

\begin{thm}\label{thm-F}
    Except for $(t,n)=(2,5)$, $(2,8)$, $(2,11)$ and $(2,14)$, the coefficient of $q^n$ in $F$ is non-negative.
\end{thm}

\pf Clearly, for $t\geq3$, the coefficient of $q^n$ in $F$ is non-negative. So we just need to consider the case of $t=2$.

For $t=2$, we have
\begin{equation}
\begin{aligned}
    F=&\frac{(q^2;q^2)_\infty}{(q;q)_\infty}\frac{q^{3}(1+q^3)(1-q)(1-q^2)}{1-q^{6}}\\
    =&\frac{(q^2;q^2)_\infty}{(q;q)_\infty}\frac{(1-q)(q^2-q^3)}{1-q^{6}}q(1+q)(1+q^3)
\end{aligned}
\end{equation}
Using \eqref{equ-G(q)}, we have that the coefficient of $q^n$ in $\dfrac{(q^2;q^2)_\infty}{(q;q)_\infty}\dfrac{(1-q)(q^2-q^3)}{1-q^{6}}$ is non-negative coefficients for $n\geq 22$. Thus for $n\geq 27$, the coefficient of $q^n$ in $F$ is non-negative. After checking $n\leq26$, we have that $F$ has non-negative coefficients of $q^n$ except for $(t,n)=(2,5)$, $(2,8)$, $(2,11)$ and $(2,14)$.\qed

{\noindent \it Proof of Theorem \ref{thm-bt2-bt3}.} From Theorem \ref{thm-D}, Theorem \ref{thm-E} and Theorem \ref{thm-F}, we see that  $b_{t,2}(n)\ge b_{t,3}(n)$ for $t\geq3$ and $n\geq 4$. And when $t=2$, again by Theorem \ref{thm-D}, Theorem \ref{thm-E} and Theorem \ref{thm-F}, we deduce that $b_{2,2}(n)\ge b_{2,3}(n)$ for $n\ge 15$.  It is routine to check that for $0\le n\le 14$,
$b_{2,2}(n)\ge b_{2,3}(n)$ still holds. This proves Theorem \ref{thm-bt2-bt3}.\qed

\end{document}